\documentclass[titlepage,11pt]{article}
\usepackage{latexsym,tikz}
\usetikzlibrary{graphs}
\usepackage{amsfonts}
\usepackage{amsmath}
\newcommand\blackslug{\hbox{\hskip 1pt \vrule width 4pt height 8pt depth 1.5pt
        \hskip 1pt}}
\newcommand\bbox{\hfill \quad \blackslug \bigbreak}
\def\d{\hbox{-}}
\def\c{\hbox{-}\cdots\hbox{-}}
\def\ll{,\ldots,}

\title{Induced subgraphs of graphs with large chromatic number.\\
XII. Distant stars}
\author{Maria Chudnovsky\thanks{Supported by NSF grant DMS-1550991 and
US Army
Research Office Grant W911NF-16-1-0404.}\\
Princeton University, Princeton, NJ 08544
\\
\\
Alex Scott\thanks{Supported by a Leverhulme Trust Research Fellowship.}\\
Oxford University, Oxford, UK
\\
\\
Paul Seymour\thanks{Supported by ONR grant N00014-14-1-0084 and 
NSF grant DMS-1265563.}\\
Princeton University, Princeton, NJ 08544}

\date{}

\newtheorem{thm}{}[section]

\newcommand{\Proof}{\noindent{\bf Proof.}\ \ }

\begin{document}
\maketitle
\begin{abstract}
The Gy\'arf\'as-Sumner conjecture asserts that if $H$ is a tree then every graph with bounded clique number and 
very large chromatic number contains $H$ as an induced subgraph. This is still open, although it has been proved for a few simple families of trees,
including trees of radius two, some special trees of radius three, and subdivided stars.
These trees all have the property that their vertices of degree more than two are clustered quite closely together.  In this paper,
we prove the conjecture for two families of trees which do not have this restriction.  
As special cases, these families contain all double-ended brooms and two-legged caterpillars.

%
%

\end{abstract}

\section{Introduction}

All graphs in this paper are finite and simple.  If $G$
is a graph, then $\chi(G)$ denotes its chromatic number, and $\omega(G)$ denotes its clique number, that is, the cardinality
of the largest clique of $G$. 

Let $H$ be a graph.  When is there a function $f$ such that $\chi(G)\le f(\omega(G))$
for every graph $G$ not containing $H$ as an induced subgraph? Let us call such a graph $H$ {\em $\chi$-bounding}.
Every $\chi$-bounding graph $H$ is a forest, because we could take
$G$ to have large girth and large chromatic number, and every such graph $G$ should contain $H$. 
The Gy\'arf\'as-Sumner conjecture~\cite{gyarfastree, sumner} asserts that the converse holds:
\begin{thm}\label{gyarfasconj}
{\bf Conjecture: } Every forest is $\chi$-bounding.
\end{thm}

It is easy to see that a forest is $\chi$-bounding if and only if all its components are $\chi$-bounding,
so the question reduces to trees.  Despite considerable attention, there are still only a few families of trees for which the Gy\'arf\'as-Sumner conjecture has been proved.
The only trees that have been shown to be $\chi$-bounding so far are:
\begin{itemize}
\item trees of radius at most two~(Gy\'arf\'as, Szemer\'edi and Tuza \cite{gst} in the triangle-free case; Kierstead and Penrice \cite{kp} in the general case);
\item trees that can be obtained from a tree of radius at most two by subdividing once {\em every} edge incident with the root~(Kierstead and Zhu \cite{kz}); and
\item subdivisions of stars (this follows from the 
``topological" version of the  Gy\'arf\'as-Sumner conjecture proved in \cite{scott}:
for every tree $T$ there is a function $f$ such that $\chi(G)\le f(\omega(G))$
for every graph $G$ containing no subdivision of $T$ as an induced subgraph.  In fact, it is enough to exclude the finite family of subdivisions of $T$ such that each edge is subdivided at most $c_T$ times, where $c_T$ is a constant depending only on the radius of $T$).  
\end{itemize}
In addition, two of us hope to show in a later paper~\cite{newbrooms} that 
every tree is $\chi$-bounding that can be obtained from a tree of radius at most two by subdividing once {\em some of} the edges incident with the root, thus
unifying the first two classes above; but the proof of that is long and difficult. 

All the trees mentioned so far have the property that their vertices of degree greater than two are all clustered closely together.
However, the conjecture is not known for any tree that contains a distant pair of vertices with degree more than two.  The aim of this paper is to show the existence of such trees.

We begin with two special cases.
Take a six-vertex path, and for each of its two middle vertices $v$ say, add another vertex adjacent to $v$. 
We obtain a tree with eight vertices, and it was not previously known whether this tree is $\chi$-bounding. More 
generally let us say a {\em two-legged caterpillar} is a tree obtained
from a path by adding two more vertices, each with one neighbour in the path.
We will prove:
\begin{thm}\label{caterpillars}
Every two-legged caterpillar is $\chi$-bounding.
\end{thm}

\begin{figure}
\centering

\begin{tikzpicture}[scale=.8,auto=left]
\tikzstyle{every node}=[circle,draw]
\node (a0) at (0,0) {};
\node (a) at (1,0) {};
\node (b) at ( 2,0) {};
\node (c) at ( 3,0) {};
\node (d) at ( 4, 0) {};
\node (e) at ( 5, 0) {};
\node (f) at (6, 0) {};
\node (g) at (2, -1) {};
\node (h) at (4,-1) {};
\node (f1) at (7,0) {};

\node (aa) at (11,0) {};
\node (bb) at ( 12,0) {};
\node (cc) at ( 13,0) {};
\node (dd) at ( 14, 0) {};
\node (ee) at ( 15, 0) {};
\node (ff) at (16, 0) {};
\node (gg) at (11, 1) {};
\node (hh) at (11,-1) {};
\node (ii) at (16, 1) {};
\node (jj) at (16,-1) {};

\foreach \from/\to in {a0/a,a/b,b/c,c/d,d/e,e/f,f/f1,b/g,d/h,aa/bb,bb/cc,cc/dd,dd/ee,ee/ff,gg/bb,hh/bb,ii/ee,jj/ee}
\draw [-] (\from) -- (\to);
\end{tikzpicture}

\caption{A 2-legged caterpillar and a double broom} \label{fig:1}
\end{figure}
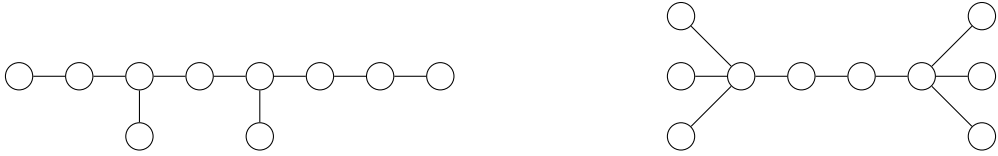

A {\em star} is a tree in which one vertex is adjacent to all the others, and  a
{\em broom} is a tree obtained from a star by replacing one of its edges by a path of arbitrary length.
A tree is a {\em subdivided star} if it has at most one vertex of degree at least three.
(All brooms are $\chi$-bounding~\cite{gyarfas}, and indeed they are subdivided stars.)
A {\em double broom} is a tree obtained from two disjoint stars by adding a path between the 
centres of the stars. We will show:
\begin{thm}\label{doublebroom} 
Every double broom is $\chi$-bounding.
\end{thm}

Our two main theorems say that certain types of trees are $\chi$-bounding, and the two results above are special cases.
The first, implying \ref{caterpillars}, is
\begin{thm}\label{maintwig}
Let $H$ be a tree obtained by adding one vertex to a subdivided star. Then $H$ is $\chi$-bounding.
\end{thm}

The second, implying \ref{doublebroom}, is:
\begin{thm}\label{mainbroom}
Let $H$ be a tree obtained from a subdivided star and a star by adding a path joining their centres. Then $H$ is $\chi$-bounding.
\end{thm}

The proofs of the two results are almost the same, and we will prove them together. 
The trees of \ref{maintwig} and \ref{mainbroom} both have the property that they have only two vertices with
degree more than two, and the length  ($k$ say) of the path between these two vertices turns out to be the key parameter.
The proof method depends
on whether there is a ball of radius at most $k$ and large chromatic number or not, so let us make that precise.
If $v$ is a vertex of a graph $G$, $N^{k}(v)$ or $N^{k}_G(v)$ denotes the set of vertices of $G$ with distance exactly $k$
from $v$, and $N^{k}[v]$ or $N^{k}_G[v]$ denotes the set with distance at most $k$ from $v$.
We sometimes write $\chi(X)$ for $\chi(G[X])$ when there is no risk of ambiguity.
If $G$ is a nonnull graph  and $k\ge 1$,
we define $\chi^{k}(G)$ to be the maximum of $\chi(N^{k}[v])$ taken over all vertices $v$ of $G$.
(For the null graph $G$ we define $\chi^{k}(G)=0$.)
If $H$ is a subgraph of $G$ and $u,v\in V(H)$, the distance between $u,v$ in $H$ may be greater than the distance between $u,v$ in $G$,
and the {\em $H$-distance between $u,v$} means the distance between $u,v$ in $H$.

Let $d\ge 1$, and take a $d$-star (that is, a copy of the complete bipartite graph $K_{1,d}$). Now subdivide 
each of its edges $d-1$ times; that is, replace each edge by a path of length $d$, joining the same pair of vertices, 
and internally pairwise disjoint. This produces a subdivided star and we call it a {\em $d$-superstar}. 
Let $k\ge 1$; we define a ``$(k,d)$-binary star'' and a ``$(k,d)$-bristled star'' as follows.
\begin{itemize}
\item Take the
disjoint union of a $d$-superstar and a $d$-star, and join their centres with a path of length $k$. We call this tree
a {\em $(k,d)$-binary star}. 
\item Take the disjoint union of a $d$-superstar $S$ and a path $T$ of length $d+1$,
and join the centre of $S$ and the second vertex of $T$ with a path of length $k$. We call this tree a {\em $(k,d)$-bristled star}.
\end{itemize}
Every tree $H$ as in \ref{mainbroom} is an induced subgraph of a $(k,d)$-binary star for some $k,d\ge 1$, and every tree $H$ as in 
\ref{maintwig} is an induced subgraph of a $(k,d)$-bristled star for some $k,d\ge 1$. It therefore suffices to prove
\ref{mainbroom} and \ref{maintwig} for trees $H$ that are $(k,d)$-binary stars, and for those that are $(k,d)$-bristled stars.
Let us say a graph is {\em $(k,d)$-starry} if it has an induced subgraph that is a $(k,d)$-binary star and one that is a 
$(k,d)$-bristled star. The following is our main result, implying both \ref{maintwig} and \ref{mainbroom}.

\begin{figure}
\centering

\begin{tikzpicture}[scale=.8,auto=left]
\tikzstyle{every node}=[circle,draw]
\node (a) at (11,0) {};
\node (b) at ( 12,0) {};
\node (c) at ( 13,0) {};
\node (d) at ( 14, 0) {};
\node (e) at ( 15, 0) {};
\node (g) at (11, 1) {};
\node (h) at (11,-1) {};
\node (i) at (15, 1) {};
\node (j) at (15,-1) {};
\node (k) at (10, 1) {};
\node (l) at (10,-1) {};
\node(m) at (10,0) {};
\node (k2) at (9, 1) {};
\node (l2) at (9,-1) {};
\node(m2) at (9,0) {};

\node (aa) at (1,0) {};
\node (bb) at ( 2,0) {};
\node (cc) at ( 3,0) {};
\node (dd) at ( 4, 0) {};
\node (ee) at ( 5, 0) {};
\node (ff) at (6, 0) {};
\node (gg) at (1, 1) {};
\node (hh) at (1,-1) {};
\node (ii) at (4, 1) {};
\node (jj) at (7,0) {};
\node (kk) at (0, 1) {};
\node (ll) at (0,-1) {};
\node(mm) at (0,0) {};
\node (kk2) at (-1, 1) {};
\node (ll2) at (-1,-1) {};
\node(mm2) at (-1,0) {};

\foreach \from/\to in {a/b,b/c,c/d,d/e,g/b,h/b,i/d,j/d,k/g,l/h,m/a,aa/bb,bb/cc,cc/dd,dd/ee,ee/ff,gg/bb,hh/bb,ii/dd,jj/ff,kk/gg,ll/hh,mm/aa,k/k2,l/l2,m/m2,
kk/kk2,ll/ll2,mm/mm2}
\draw [-] (\from) -- (\to);
\end{tikzpicture}

\caption{A $(3,2)$-bristled star and a $(3,2)$-binary star} \label{fig:2}
\end{figure}
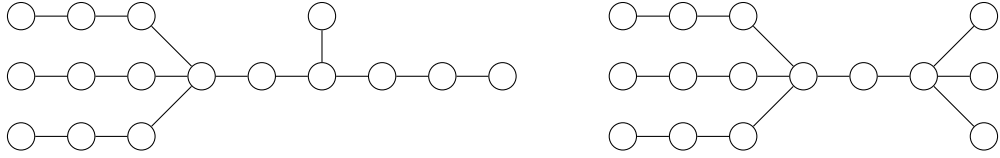

\begin{thm}\label{mainthm}
For all $\kappa\ge 0$ and $k,d\ge 1$, there exists $c\ge 0$ such that every graph $G$ with $\omega(G)\le \kappa$
and $\chi(G)>c$ is $(k,d)$-starry.
\end{thm}

The proof of \ref{mainthm} is given at the end of section six.

\section{Using criticality}

Our main tool is a set of lemmas proved in this section, that if $X$ is a subset of $V(G)$ of small chromatic number,
and $G$ itself has large chromatic number, and deleting $X$ from $G$ reduces the chromatic number, then there are useful
subgraphs rooted at some vertex in $X$ and growing out into $G\setminus X$. We begin with:

\begin{thm}\label{creatures}
Let $d\ge 0$ be an integer, let $G$ be a graph with chromatic number more than $d$, and let $X\subseteq V(G)$ be stable,
such that $\chi(G\setminus X)<\chi(G)$. Then some vertex in $X$ has at least $d$ neighbours in $V(G)\setminus X$. 
\end{thm}
\Proof
Let $\chi(G)=k+1$, and so $k\ge d$. Let $\phi:V(G)\setminus X\rightarrow \{1\ll k\}$ be a $k$-colouring of $G\setminus X$.
For each $x\in X$,
if $x$ has at most $k-1$ neighbours in $V(G)\setminus X$ then we may choose $\phi(x)\in \{1\ll k\}$, different from $\phi(v)$
for each neighbour $v\in V(G)\setminus X$ of $x$; and this extends $\phi$ to a $k$-colouring of $G$, which is impossible. Thus
for some $x\in X$, $x$ has at least $k\ge d$ neighbours in $V(G)\setminus X$. 
This proves \ref{creatures}.~\bbox

If $X\subseteq V(G)$, let us say an {\em $X$-split} is a triple $(x,y,Z)$, where 
\begin{itemize}
\item $x\in X$, $y\in V(G)\setminus X$, and $Z\subseteq V(G)\setminus (X\cup \{y\})$;
\item $x$ is adjacent to $y$ and has at least one neighbour in $Z$;
\item $y$ has no neighbours in $Z$; and
\item $G[Z]$ is connected.
\end{itemize}
Let us say $\chi(Z)$ is the {\em chromatic number} of an $X$-split $(x,y,Z)$.
Next we need:

\begin{thm}\label{stabletwiglemma}
For all $c,\tau\ge 0$ there exists $c'$ with the following property. Let $G$ be a graph with chromatic number more than $c'$,        
such that $\chi^1(G)\le \tau$;
and let $X\subseteq V(G)$ be stable, such that $\chi(G\setminus X)<\chi(G)$. Then there is an $X$-split
in $G$ with chromatic number more than $c$.
\end{thm}
\Proof
Let $c'=(2c+3\tau+3)\tau$, let $G$ be a graph with chromatic number more than $c'$ and $\chi^1(G)\le \tau$, 
and let $X\subseteq V(G)$ be stable, such that 
$\chi(G\setminus X)<\chi(G)$. We prove the result by induction on $|X|$. Choose $x\in X$. If $\chi(G\setminus x)=\chi(G)$,
let $G'=G\setminus x$, and $X'=X\setminus \{x\}$; then $\chi(G'\setminus X')<\chi(G')$, and the result follows from the
inductive hypothesis. Thus we may assume that $\chi(G\setminus x)<\chi(G)$. 

Let $k=\chi(G)-1$, and let $\phi:V(G)\setminus \{x\}\rightarrow\{1\ll k\}$ 
be a $k$-colouring of $G\setminus x$.
Let $N$ denote the set of all neighbours of $x$.
Then, since $G$ does not admit a $k$-colouring, it follows that for $1\le i\le k$ there exists $n\in N$ with $\phi(n)=i$.
Let $\overline{G}$ be the complement graph of $G$, and
let $N_1\ll N_t$ be the vertex sets of the components of $\overline{G}[N]$. 
Now $\chi(N)\le \tau$, since $\chi^1(G)\le \tau$.
But for $1\le i<j\le t$, every vertex of $N_i$ is adjacent to every vertex of $N_j$, and so $\chi(N)=\sum_{1\le i\le t}\chi(N_i)$.
Consequently $\sum_{1\le i\le t}\chi(N_i)\le \tau$, and in particular, $t\le \tau$ since each $\chi(N_i)>0$. 
For $1\le i\le t$ let $D_i$ be the set $\{\phi(v):v\in N_i\}$, that is, the set of colours that appear in $N_i$. Now
$D_1\ll D_n$ are pairwise disjoint and have union $\{1\ll k\}$, so we may assume that $|D_1|\ge k/t\ge 2c+3\tau+3$.

Let $Y$ be the set of vertices of $G\setminus N_1$ that are adjacent to every vertex in $N_1$; thus $x\in Y$, and
$N\setminus N_1\subseteq Y$. Let $W$ be the set of vertices of $G\setminus (N_1\cup Y)$ that have neighbours in $N_1$.
\\
\\
(1) {\em If $\chi(W)>c+2\tau+1$ then the theorem holds.}
\\
\\
We assume that $\chi(W)>c+2\tau+1$. Since $X$ is stable, it follows that $\chi(W\setminus X)>c+2\tau$; and since $\chi(N)\le \tau$,
$\chi(W\setminus (N\cup X))>c+\tau$. Choose $q\in N_1$; then since $\chi(N^1(q))\le \tau$, the set of vertices in $W\setminus (N\cup X)$
that are nonadjacent to $q$ has chromatic number more than $c$, and so there exists $Z\subseteq W\setminus (N\cup X)$
with $\chi(Z)>c$, such that $G[Z]$ is connected and $q$ has no neighbour in $Z$. Let $P$ be the set of vertices in $N_1$
that have a neighbour in $Z$, and $Q=N_1\setminus P$; then $P\ne\emptyset$, since every vertex of $Z$ has a neighbour in $N_1$,
and $Q\ne \emptyset$, since $q\in Q$. Since $\overline{G}[N_1]$ is connected, there exist $y,z\in N_1$, nonadjacent,
such that $z\in P$ and $y\in Q$. But then $(x,y,Z\cup \{z\})$ is an $X$-split satisfying the theorem.
This proves (1).

\bigskip

In view of (1), we assume henceforth that $\chi(W)\le c+2\tau+1$.
\\
\\
(2) {\em Let $C$ be the vertex set of a component of $G\setminus (N_1\cup Y)$. If $C\cap W\ne \emptyset$ and
$\chi(C)>2c+2\tau+2 $ then the theorem holds.}
\\
\\
For then $\chi(C\setminus (W\cup X))>c$, and so there is a subset $Z\subseteq C\setminus (W\cup X)$
with $\chi(Z)>c$ such that $G[Z]$ is connected. Since $Z\subseteq C$ and $C\cap W\ne \emptyset$, there is
a path of $G[C]$ between $W$ and $Z$; choose such a path, $P$ say, minimal. Now $P$ has length at least one,
since $Z\cap W=\emptyset$. Let $w$ be the end of $P$ in $W$. It follows that no vertex of $P$ different from $w$
has a neighbour in $N_1$. Since $\overline{G}[N_1]$ is connected and
$w$ has a neighbour and a non-neighbour in $N_1$, there exist nonadjacent 
$y,z\in N_1$ such that $w$ is adjacent to $z$ and not to $y$.

Suppose first that no vertex of $P$ belongs to $X$. 
Then $(x,y,Z\cup V(P)\cup\{z\})$ is the desired $X$-split. We may assume therefore that some vertex
of $P$ belongs to $X$. Choose $x'\in X\cap V(P)$ such that the subpath of $P$ ($P'$ say) between $x'$ and $Z$ is minimal. 
If $x'\ne w$, let $y'$ be the vertex of $P$ adjacent to $x'$ that does not belong to $V(P')$,
and if $x'=w$, let $y'=z$. In either case, $y'\notin X$, since $X$ is stable. 
Then $\{x',y',(V(P')\setminus \{x'\})\cup Z)$ is the desired $X$-split. This proves (2).

\bigskip

Let $U$ be the union of the vertex sets of all components of $G\setminus (N_1\cup Y)$ that have nonempty intersection with $W$.
\\
\\
(3) {\em If $v\in V(G)\setminus (N_1\cup U\cup \{x\})$
and has a neighbour $u\in N_1\cup U\cup \{x\}$, then $\phi(v)\notin D_1$.}
\\
\\
We claim that $v\in Y$; for if $u\in N_1$ then $v$ has a neighbour in $N_1$ and $v\notin U\cup N_1\cup \{x\}$, so $v\in Y$; if
$u\in U$ then $v\in N_1\cup Y$ since $v\notin U$, and so again $v\in Y$; and if $u=x$ then $v\in N^1(x)\subseteq N_1\cup Y$,
and again $v\in Y$. Thus $v\in Y$. Consequently $\phi(v)$ is different from $\phi(w)$
for every $w\in N_1$; and so from the definition of $D_1$ it follows that $\phi(v)\notin D_1$. 
This proves (3).

\bigskip
By (2) we may assume (for a contradiction) that $\chi(U)\le 2c+2\tau+2$. It follows that $\chi(N_1\cup U\cup \{x\})\le 2c+3\tau+3$.
Consequently $G[N_1\cup U\cup \{x\}]$ admits a colouring $\psi$ using only the colours in $D_1$,
since $|D_1|\ge 2c+3\tau+3$. For each $v\in V(G)$, define $\phi'(v) = \phi(v)$ if $v\notin N_1\cup U\cup \{x\}$, and $\phi'(v)=\psi(v)$
if $v\in N_1\cup U\cup \{x\}$. 
By (3) this gives a $k$-colouring of $G$, which is impossible.
This proves \ref{stabletwiglemma}.~\bbox

To use this, we combine it with a version of Gy\'arf\'as' path theorem (see \cite{longoddholes} for this version):
\begin{thm}\label{getpath}
Let $G$ be a graph, let $k\ge 0$, let $C\subseteq V(G)$, and let $x_0\in V(G)\setminus C$,
such that
$G[C]$ is connected,
$x_0$ has a neighbour in $C$, and
$\chi(C)>k \chi^1(G)$.
Then there is an induced path $x_0 \c x_{k}$ of $G$ where $x_1\ll x_{k}\in C$, and a subset $C'$ of $C$,
with the following properties:
\begin{itemize}
\item $x_0\ll x_{k}\notin C'$;
\item $G[C']$ is connected;
\item $x_{k}$ has a neighbour in $C'$, and $x_0\ll x_{k-1}$ have no neighbours in $C'$; and
\item $\chi(C')\ge \chi(C)-k\chi^1(G)$.
\end{itemize}
\end{thm}

We deduce:
\begin{thm}\label{twiglemma}
For all $d, \tau\ge 0$ there exists $c'$ with the following property. Let $G$ be a graph with chromatic number more than $c'$,
such that $\chi^1(G)\le \tau$; and
let $X\subseteq V(G)$ be stable, such that $\chi(G\setminus X)<\chi(G)$. Then there is a vertex $x\in X$,
an induced path $P$ of length $d$ with one end $x$ and no other vertices in $X$, and a vertex $y\in V(G)\setminus X$ that is
adjacent to $x$ and has no other neighbour in $V(P)$.
\end{thm}
\Proof Let $c=d\tau$ and let $c'$ satisfy \ref{stabletwiglemma}. We claim that $c'$ satisfies \ref{twiglemma}. For let
$G, X$ be as in the theorem, with $\chi(G)>c'$. By \ref{stabletwiglemma} there is an $X$-split $(x,y,Z)$ with chromatic number 
more than $c=d\tau$. By \ref{getpath} there is an induced path $P$ of $G[Z\cup \{x\}]$ with one end $x$, of length $d$.
But then $x,y,P$ satisfy \ref{twiglemma}. This proves \ref{twiglemma}.~\bbox

We combine \ref{creatures} and \ref{twiglemma} in the following. If $X\subseteq V(G)$, we say that $x\in X$ is
{\em $d$-equipped in $V(G)\setminus X$} if $x$
has at least $d$ neighbours in $V(G)\setminus X$, pairwise nonadjacent, and there is 
an induced path $P$ of length $d$ with one end $x$ and no other vertices in $X$, and a vertex $y\in V(G)\setminus (X\cup V(P))$ 
that is
adjacent to $x$ and has no other neighbour in $V(P)$. We speak of the set of $d$ pairwise nonadjacent neighbours,
the path $P$ and the vertex $y$ as the {\em parts of the equipment} of $x$.

\begin{thm}\label{bigstabletwiglemma}
For all $d, \tau\ge 0$ there exists $c$ with the following property. Let $G$ be a graph with chromatic number more than $c$,
such that $\chi^1(G)\le \tau$; and
let $X\subseteq V(G)$ be stable, such that $\chi(G\setminus X)<\chi(G)$. Then there is a vertex $x\in X$
that is $d$-equipped in $V(G)\setminus X$.
\end{thm}
\Proof
Choose $d'\ge 0$ such that every graph with $d'$ vertices has either a clique of cardinality $\tau+1$ or a 
stable set of size $d$. Choose $c\ge d'$ such that \ref{twiglemma} holds with $c'$ replaced by $c$. 
Let $X_1$ be the set of vertices in $X$ that have at least $d'$ neighbours in $V(G)\setminus X$; and let $X_2$
be the set of vertices $x\in X$ such that there is an induced path $P$ of length $d$ with one end $x$ and 
no other vertices in $X$, and a vertex $y\in V(G)\setminus X$ that is
adjacent to $x$ and has no other neighbour in $V(P)$. Suppose for a contradiction that $X_1\cap X_2=\emptyset$.
Since $X$ is stable, it follows that no vertex in $X_2$ has $d'$ neighbours in $V(G)\setminus X_2$; and so by
\ref{creatures}, it follows that $\chi(G)=\chi(G\setminus X_2)$. By \ref{twiglemma} applied to $X\setminus X_2$ in the graph
$G\setminus X_2$, it follows that 
$$\chi(G\setminus X_2) = \chi((G\setminus X_2)\setminus (X\setminus X_2)),$$
contradicting that $\chi(G\setminus X)<\chi(G)$. Hence there exists $x\in X_1\cap X_2$.
The choice of $d'$ implies that $x$ has $d$
pairwise nonadjacent neighbours in $V(G)\setminus X$, so 
$x$ is $d$-equipped in $V(G)\setminus X$. This proves \ref{bigstabletwiglemma}.~\bbox

Next we need a version in which $X$ may not be stable:
\begin{thm}\label{bigtwiglemma}
For all $a,d, \tau\ge 0$ there exists $c$ with the following property. Let $G$ be a graph with chromatic number more than $c$,
such that $\chi^1(G)\le \tau$; and
let $X\subseteq V(G)$ with $\chi(X)\le a$, such that $\chi(G\setminus X)<\chi(G)$. Then there is a vertex $x\in X$
that is $d$-equipped in $V(G)\setminus X$.
\end{thm}
\Proof
Choose $b$ such that \ref{bigstabletwiglemma}  holds with $c$ replaced by $b$; and let $c=ab$. We claim that $c$ 
satisfies the theorem.
Let $G$ be a graph with chromatic number more than $c$,
and let $X\subseteq V(G)$ with $\chi(X)\le a$, such that $\chi(G\setminus X)<\chi(G)$. Let
$\chi(G)=k+1$, and let $\phi:V(G)\setminus X\rightarrow \{1\ll k\}$ be a $k$-colouring of $G\setminus X$. Let $(X_1\ll X_a)$ be a partition of $X$ into $a$ stable sets.
For $1\le i\le a$ let $A_i = \{(i-1)b+1\ll ib\}$, and let $Y_i$ be the set of vertices $v\in V(G)\setminus X$
with $\phi(v)\in A_i$. If $\chi(X_i\cup Y_i)\le b$ for each $i$ then $\chi(G)=\chi(G\setminus X)$, a contradiction;
so we may assume that $\chi(X_1\cup Y_1)>b$. Since $\chi(Y_1)\le b$, and $X_1$ is stable, the choice of $b$ implies
that there is a vertex $x\in X_1$
that is $d$-equipped in $Y_1\setminus X_1$ and hence in $V(G)\setminus X$.
This proves \ref{bigtwiglemma}.~\bbox

We can iterate this (and we also throw in a bounded set of ``forbidden vertices'' $B$, but for most of 
the applications $B = \emptyset$).

\begin{thm}\label{doublecreatures}
For all $a,b,d, \tau\ge 0$ there exists $c$ with the following property. Let $G$ be a graph with chromatic number more than $c$,
such that $\chi^1(G)\le \tau$; and
let $X\subseteq V(G)$ such that $\chi(G\setminus X)<\chi(G)$ and $\chi(X)\le a$. 
Let $B\subseteq V(G)\setminus X$ with $|B|\le b$.
Let $M$ denote the set of vertices in $V(G)\setminus X$ with a neighbour in $X$.
Then there is a vertex $x\in X$
that is $d$-equipped in $V(G)\setminus X$, and there is a neighbour $x'$ of $x$ in $M\setminus B$,
such that $x'$ is $d$-equipped in $B\cup (V(G)\setminus (X\cup N^1(x)))$.
\end{thm}
\Proof By increasing $d$, we may assume that $d>b$.
Choose $c_1$ such that \ref{bigtwiglemma} holds with $a,c$ replaced by $\tau,c_1$.
Choose $c$ such that \ref{bigtwiglemma} holds with $a,c$ replaced by $c_1+a+b,c$.
Now let $G,X,M$ be as in the theorem with $\chi(G)>c$. 
Let $X_1$ be the set of vertices in $X$ that are
$d$-equipped in $V(G)\setminus X$, and let $M_1$
be the set of vertices in $V(G)\setminus X$ with a neighbour in $X_1$. 

Suppose first that $\chi(X\cup M_1)\le c_1+a+b$. 
Let $X'= (X\cup M_1)\setminus B$; then $X\subseteq X'$, and so $\chi(G\setminus X')<\chi(G)$. From
\ref{bigtwiglemma}, some vertex $x'\in X'$ is $d$-equipped in $V(G)\setminus X'$.
Now
$x'\notin X\setminus X_1$, since no vertex in $X\setminus X_1$ is $d$-equipped in $V(G)\setminus X$.
Also $x'$ has at least $d$ neighbours in $(V(G)\setminus (X\cup M_1))\cup B$, and therefore at least one
neighbour in $V(G)\setminus (X\cup M_1)$, since $|B|\le b<d$; and no vertex in $X_1$ has any neighbour in 
$V(G)\setminus (X\cup M_1)$. Consequently $x'\notin X_1$, and so $x'\in M_1\setminus B$. Choose $x\in X_1$ adjacent to $x'$.
Then $N^1(x)\setminus X\subseteq M_1$, and so $x'$ is $d$-equipped in $B\cup (V(G)\setminus (X\cup N^1(x)))$.

We may assume therefore that $\chi(X\cup M_1)> c_1+a+b$, and so $\chi(M_1\setminus B)>c_1$.
Choose $Z\subseteq M_1\setminus B$ minimal with $\chi(Z)>c_1$. Choose $x\in X_1$ with a neighbour in $Z$, and let $Y$ be the set
of neighbours of $x$ in $Z$. Then $\chi(Y)\le \tau$, and so by \ref{bigtwiglemma} applied to $G[Z], Y$, there exists
$x'\in Y$ that is $d$-equipped in $Z\setminus Y$. Since $x$ has no neighbours in $Z\setminus Y$, it follows that
$x,x'$ satisfy the theorem. This proves \ref{doublecreatures}.~\bbox

\section{$k$-balls with large chromatic number}

Let $k,d\ge 1$ be integers. A {\em $(k,d)$-broom } is a tree obtained from a path $v_0\d v_1\c v_k$ by adding $d$ new vertices, 
each adjacent to $v_k$.  We call $v_0$ the {\em root} of the broom.
A {\em $(k,d)$-bristle} is obtained from a path $v_0\d v_1\c v_{k+d}$ by adding one new vertex adjacent to $v_k$. We call $v_0$
the {\em root} of the bristle. By a {\em $(k,d)$-broom in $G$} we mean a $(k,d)$-broom that
is an induced subgraph of $G$, and we use similar language for other kinds of tree.

\begin{thm}\label{firstbroom}
Let $k,d,\tau\ge 1$; then there exists $c$ with the following property.
Let $G$ be a graph with $\chi^1(G)\le \tau$, and let $z\in V(G)$, such that $\chi(N^k[z])>c$. 
Then there is a $(k,d)$-broom and a $(k,d)$-bristle in $G$, both with root $z$.
\end{thm}
\Proof
Let $c_1=2\tau$, and inductively for $i\ge 2$, choose $c_i\ge 2dc_{i-1}$ such that 
\ref{bigtwiglemma} is satisfied
with $a,c$ replaced by $c_{i-1}, c_i/2$ respectively.
We prove by induction on $k$ (for the given value of $d$) that setting $c=c_k$
satisfies the theorem. Thus we may assume that either $k=1$ or the claim holds for $k-1$.

Let $G$ satisfy $\chi^1(G)\le \tau$, and let $z\in V(G)$, such that $\chi(N^k[z])>c_k$. For each $s\ge 0$, let $L_s=N^s(z)$.
Since $\chi(N^k[z])>c_k$, there exists $s\ge 0$ with $s\le k$ such that $\chi(L_s)>c_k/2$.
Since $\chi(L_0),\chi(L_1)\le \tau$, it follows that $s\ge 2$; and so $k\ge 2$, and the claim holds for $k-1$.

Choose $S\subseteq L_s$ minimal such that $\chi(S)=\chi(L_s)$. 
Choose $u\in S$, and choose $v\in L_1$ joined to $u$ by a path of length $s-1$.

Let $G'$ be the graph $G\setminus (L_0\cup L_1\setminus \{v\})$. If $\chi(N^{k-1}_{G'}[v])>c_{k-1}$, 
then from the inductive
hypothesis applied to $G'$, there is a $(k-1,d)$-broom and a $(k-1,d)$-bristle in $G'$ with root $v$.
But then adding the edge $zv$ gives the desired $(k,d)$-broom  and $(k,d)$-bristle in $G$ with root $z$. 
We may therefore assume that $\chi(N^{k-1}_{G'}[v])\le c_{k-1}$.

Let $X$ be the set of vertices in $S$ that have $G'$-distance at most $k-1$ from $v$.
Thus $X\subseteq N^{k-1}_{G'}[v]$, and so $\chi(X)\le c_{k-1}$.
Now $u\in X$, since $s\le k$, and so 
$X\ne \emptyset$. From the minimality of $S$, $\chi(S\setminus X)<\chi(S)$.

By \ref{bigtwiglemma}, since $\chi(S)\ge c_k/2\ge c_{k-1}d$, and $\chi(X)\le c_{k-1}$, it follows that some vertex $x\in X$
is $d$-equipped in $S\setminus X$. 
Since $x\in X$, there is an induced
path $P$ of length at most $k-1$ such that
$V(P)\cap (L_0\cup L_1)=\{v\}$. Since $x$ has a neighbour in $S\setminus X$, it follows that the length of $P$ is 
exactly $k-1$, and no vertex of $P$ different from $x$ has a neighbour in $S\setminus X$. Also $z$
has no neighbours in $S\setminus X$ since $s\ge 2$.
But then $P$ together with the edges $zv$ and the various parts of the equipment of $x$ gives a $(k,d)$-broom 
and a $(k,d)$-bristle, both with root $z$.
This proves \ref{firstbroom}.~\bbox

If $A,B$ are disjoint subsets of $V(G)$, we say that $A$ {\em covers} $B$ if every
vertex in $B$ has a neighbour in $A$.

\begin{thm}\label{treesplit}
For all $c,\tau\ge 0$ and $d\ge 1$ and $k\ge r\ge 2$, there exists $c'\ge 0$ with the following property. Let $G$ be a graph
with $\chi^{r-1}(G)\le \tau$, and let 
$z\in V(G)$ such that 
$\chi(N^r_G(z))>c'$. Then there is a $(k,d)$-broom $H_1$
and a $(k,d)$-bristle $H_2$ in $G$, both with root $z$, and a subset $W$ of $V(G)$, such that 
\begin{itemize}
\item $z\in W$, and no other vertex of $H_1\cup H_2$ belongs to $W$ or has a neighbour in $W\setminus \{z\}$;
\item $\chi(N^r_{G'}(z))>c$, where $G'$ is the graph $G[W]$.
\end{itemize}
\end{thm}
\Proof
Choose $c_1$ such that \ref{firstbroom} is satisfied with $c$, $k$ replaced by $c_1, k$.
Let
$$c'=((2k+2d+3)\tau+c+c_1)\tau.$$
Now let $G,z$ be as in the theorem.
For $u\in N^1(z)$ and $v\in N^i(z)$ where $i\ge 1$, we say that $u$ is an {\em ancestor} of $v$ and $v$ is a {\em descendant}
of $u$ if there is a path of length $i-1$
between $u,v$.
Since $\chi(N^1(z))\le \tau$, there is a partition of $N^1(z)$ into $\tau$ stable sets; and since every vertex in $N^r(z)$
has an ancestor in $N^1(z)$, there is a stable set $L_1\subseteq N^1(z)$ such that, if $L_r$ denotes the set of 
vertices in $N^r(z)$ with an ancestor in $L_1$, then
$$\chi(L_r)>c'/\tau= (2k+2d+3)\tau+c+c_1.$$
For $2\le i\le r-1$, let $L_i$ be the set of descendants in $N^i(z)$ of members of $L_1$.

Choose $B\subseteq L_1$ maximal such that the set of vertices in $L_r$ with no ancestor in $B$
has chromatic number at least $c_1$. 
Let $L_0'=\{z\}$ and for $1\le i\le r$, let $L_i'$ be the set of vertices in $L_i$
with no ancestor in $B$. It follows that $L_1'=L_1\setminus B$, and $\chi(L_r')>c_1$, and for $1\le i\le r$, $L_{i-1}'$
covers $L_i'$. Let 
$$G_1 = G[L_0'\cup\cdots\cup L_r'].$$ 
Then $\chi(N^r_{G_1}(z))>c_1$, and so $\chi(N^{k}_{G_1}[z])>c_1$.
It follows from \ref{firstbroom} applied to $G_1$, and from the choice of $c_1$, that 
there is a $(k,d)$-broom $H_1$ and a $(k,d)$-bristle $H_2$ in $G$, both with root $z$, and both with vertex set
a subset of $L_0'\cup\cdots\cup L_r'$.

In particular, there is a vertex in $L_1\setminus B$, $v$ say, and since the set of descendants of $v$ in $L_r$
has chromatic number at most $\tau$, the maximality of $B$ implies that $\chi(L_r')\le c_1+\tau$. 
Consequently 
$$\chi(L_r\setminus L_r')> 2(k+d+1)\tau+c.$$ 
Since 
$$|V(H_1)\cup V(H_2)|\le 2(k+d+1),$$
the set of vertices with $G$-distance at most $r-1$
from some vertex in $H_1\cup H_2$ has chromatic number at most $2(k+d+1)\tau$. Consequently there exists 
$L_r''\subseteq L_r\setminus L_r'$ with 
$\chi(L_r'')> c$
such that every vertex in $L_r''$ has $G$-distance at least $r$ from every vertex of $H_1\cup H_2$. Let $L_0''=\{z\}$, and 
for $1\le i \le r-1$,
let $L_i''$ be the set of vertices in $L_i\setminus L_i'$ that have $G$-distance at least $i$ from every vertex of $H_1\cup H_2$.
We claim that $L_{i-1}''$ covers $L_{i}''$ for $1\le i\le r$. Certainly $L_0''$ covers $L_1''$, so we may assume that $i\ge 2$.
Let $v\in L_i''$. Since $v\notin L_i'$, $v$ has an ancestor in $B$, and since $i\ge 2$ it follows
that $v$ has a neighbour $u\in L_{i-1}$ with an ancestor in $B$. Consequently $u\in L_{i-1}\setminus L_{i-1}'$.
But since the $G$-distance from $v$ to $V(H_1\cup H_2)$ is at least $i$, and $u,v$ are adjacent,
it follows that the $G$-distance from $u$ to $V(H_1\cup H_2)$ is at least $i-1$, and so $u\in L_{i-1}''$. This proves that
$L_{i-1}''$ covers $L_{i}''$ for $1\le i\le r$. Let $W=L_0''\cup L_1''\cup\cdots\cup L_r''$, and let $G'=G[W]$.
Then $\chi(N^r_{G'}(z))>c$. It remains to show that no vertex of $W\setminus \{z\}$ is adjacent to a vertex of $H_1\cup H_2$
different from $z$. Let $v\in L_i''$ say, where $i\ge 1$, and suppose $v$ is adjacent to some 
vertex $u\in V(H_1\cup H_2)\setminus\{z\}$.
Since $v$ has $G$-distance at least $i$ from $V(H_1\cup H_2)$, it follows that $i=1$, and so $v\in B$ and $u\in L_1'\cup L_2'$.
But $u\notin L_1$ since $L_1$ is stable, and $u\notin L_2'$ since no vertex of $L_2'$ 
has an ancestor in $B$,
a contradiction. Thus there is no such pair $u,v$. This proves \ref{treesplit}.~\bbox

Let $s\ge 1$, for $1\le i\le s$ let $k_i,d_i\ge 1$, and let $H_i$ be either a $(k_i,d_i)$-broom or a $(k_i,d_i)$-bristle; 
and let $H$ be the rooted tree obtained from
the disjoint union of $H_1\ll H_s$ by identifying the roots of $H_1\ll H_s$ to form the root of $H$. Let us call
such a tree $H$ the {\em rooted sum} of $H_1\ll H_s$.

\begin{thm}\label{bigball}
With notation as above, let $H$ be the rooted sum of $H_1\ll H_s$;
and let $1\le r\le \min(k_1\ll k_s)$.
For all $\tau\ge 0$ and $r\ge 1$ there exists $c$ with the following property. Let $G$ be a graph
with $\chi^1(G)\le \tau$, and 
let $z\in V(G)$ be a vertex such that $\chi(N^r[z])>c$. Then there is an induced subgraph of $G$ isomorphic to $H$ with root $z$.
\end{thm}
\Proof We may assume (by replacing $d_1\ll d_s$ by their maximum) that $d_1=\cdots=d_s=d$ say; for fixed $d$,
we proceed by induction  on $r$; and for
fixed $d,r$, by induction on $s$. If $r=1$ then the claim holds, setting $c=\tau$; so we may assume that $r>1$, and
(for the same values of $d,s$ and $k_1\ll k_s$) the result holds with $r,c$ replaced by $r-1,c_1$. 
Let $H'$ be the rooted sum of $H_1\ll H_{s-1}$ (or if $s=1$, let $H'$ be a one-vertex graph, with that vertex as root).
If $s=1$ let $c_2=0$, and if $s>1$ then by induction on $s$, we may choose $c_2$ such that 
the result holds (for the same values of $d,r$ and $k_1\ll k_{s-1}$)
with $s,H,c$ replaced by $s-1,H',c_2$. 
By \ref{treesplit} we may choose $c_3$ such that \ref{treesplit} holds with $\tau,c,k,r,c'$ replaced by
$c_1,c_2,k_s,r,c_3$. We claim that setting $c=c_1+c_3$ satisfies the theorem.

For let $G,z$ be as in the theorem, with $\chi(N^r[z])>c_1+c_3$. If there is a vertex $z'$ with $\chi(N^{r-1}[z'])>c_1$ 
then the result follows from the choice of $c_1$. Thus we may assume that $\chi^{r-1}(G)\le c_1$; and in particular
$\chi(N^r(z))>c_3$. From \ref{treesplit} and the choice of $c_3$, 
there is an induced subgraph $J$ of $G$ isomorphic to $H_s$, with root $z$,
and a subset $W$ of $V(G)$, such that
\begin{itemize}
\item $z\in W$, and no other vertex of $J$ belongs to $W$ or has a neighbour in $W\setminus \{z\}$; and
\item $\chi(N^r_{G'}(z))>c_2$, where $G'$ is the graph $G[W]$.
\end{itemize}
From the choice of $c_2$, there is an induced subgraph $J'$ of $G'$, isomorphic to $H'$, with root $z$. But then the union
of $J$ and $J'$ satisfies the theorem. This proves \ref{bigball}.~\bbox

\section{Spires and cathedrals}

Let $G$ be a graph, let $P$ be an induced path of $G$, and $A,B\subseteq V(G)$, such that:
\begin{itemize}
\item $G[A]$ is connected;
\item $A\cap B=\emptyset$;
\item $A$ covers $B$;
\item $V(P)\cap B=\emptyset$, and there is an end $z$ of $P$ in $A$ such that $V(P)\cap A=\{z\}$; and
\item no vertex in $V(P)\setminus \{z\}$ has any neighbours in $(A\cup B)\setminus \{z\}$.
\end{itemize}
In this situation we say that $\mathcal{S}=(P,A,B)$ is a {\em spire} of {\em height} $d$, where $d$ is the length of $P$.
We define $V(\mathcal{S})=A\cup B\cup V(P)$. If $C\subseteq V(G)$, we say that the spire {\em dominates} $C$
if 
\begin{itemize}
\item $C$ is disjoint from $V(\mathcal{S})$;
\item there are no edges between $A\cup V(P)$ and $C$; and 
\item $B$ covers $C$.
\end{itemize}
A {\em cathedral} is a sequence of spires $(\mathcal{S}_1\ll\mathcal{S}_n)$, such that
\begin{itemize}
\item for $1\le i<j\le n$, $V(\mathcal{S}_i)\cap V(\mathcal{S}_j)=\emptyset$; and
\item for $1\le i<j\le n$, if $u\in V(\mathcal{S}_i)$ and $v\in V(\mathcal{S}_j)$ are adjacent then $u\in B_i$ and $v\in A_j\cup B_j$
\end{itemize}
where $\mathcal{S}_i=(P_i,A_i,B_i)$ for $1\le i\le n$.
We say the cathedral is {\em free} if 
\begin{itemize}
\item for $1\le i<j\le n$, if $u\in V(\mathcal{S}_i)$ and $v\in V(\mathcal{S}_j)$ are adjacent then $u\in B_i$ and $v\in B_j$.
\end{itemize}
A cathedral has {\em height} $d$ if each of its spires has height $d$, and {\em length} $n$ if it has $n$ spires.
We say a cathedral {\em dominates} a set $C\subseteq V(G)$
if each of its spires dominates $C$. 

\begin{thm}\label{freecathedral}
For all $\tau\ge 0$ and $k,d\ge 1$ there exist $c,n\ge 0$ with the following property. 
Let $G$ be a graph with $\chi^k(G)\le \tau$, and let 
$(\mathcal{S}_1\ll\mathcal{S}_n)$ be a free cathedral in $G$ of height $d$, dominating a set $C$ with $\chi(C)>c$.
Then $G$ is $(k,d)$-starry.
\end{thm}
\Proof
Let $m=2d^2$. 
Choose $n\ge 0$ such that every graph with $n$ vertices has either a clique of cardinality $\tau+1$ or a stable set of cardinality $m$.
Choose $c_0$ such that \ref{bigtwiglemma} holds with $a,c$ replaced by $\tau,c_0$.
If $k\ge 2$, choose $c$ such that \ref{doublecreatures} holds with $a$, $b$, $c$ and $d$ replaced by $(m+1)\tau + c_0$, 
$0$, $c$ and $d$ respectively; and if $k=1$, choose $c$ such that \ref{doublecreatures} holds with $a$, $b$, $c$ and $d$ 
replaced by $(m+1)\tau + c_0$, $mn$, $c$ and $d+mn$ respectively.
Now let $G$ and $(\mathcal{S}_1\ll\mathcal{S}_n)$ be as in the theorem, and let 
$\mathcal{S}_i=(P_i,A_i,B_i)$ for $1\le i\le n$.
We may assume that $C$ is minimal such that $\chi(C)>c$, and in particular, $G[C]$ is connected.
Choose $r_0\in C$, and for $1\le i\le n$ let $b_i$ be a neighbour of $r_0$
in $B_i$. Since $P_i$ has length $d$, and $A_i$ covers $B_i$, and $G[A_i\cup V(P_i)]$ is connected, 
there is an induced path $Q_i$ of $G[A_i\cup V(P_i)\cup \{b_i\}]$ of length $d-1$
with one end $b_i$.
By the choice of $n$, we may assume that $b_1\ll b_m$ are pairwise nonadjacent.

Let $R$ be the set of vertices in $C$ that are adjacent to at least $2d$ of $b_1\ll b_m$; thus $r_0\in R$. 
Let $S$ be the set of vertices in $C$ that are adjacent to at least one and to at most $2d-1$ of $b_1\ll b_m$.
Let $L_0=R$, and let $L_1$ be the set of all vertices in $C\setminus L_0$ with a neighbour in $L_0$. Let $L_2$ be the set of all vertices
in $C\setminus (L_0\cup L_1)$ that either belong to $S$ or have a neighbour in $L_1$. For $i\ge 3$, inductively let $L_i$
be the set of vertices in $C\setminus (L_0\cup\cdots\cup L_{i-1})$ with a neighbour in $L_{i-1}$. Let $W$ be the set of vertices
in $C$ with $G$-distance at most $k$ from a vertex in $\{r_0,b_1\ll b_m\}$.
\\
\\
(1) {\em If $\chi(L_k\setminus W)>c_0$ then the theorem holds.}
\\
\\
Suppose that $\chi(L_k\setminus W)>c_0$, and choose $Y\subseteq L_k\setminus W$ minimal with $\chi(Y)>c_0$. 
In particular $(R\cup S)\cap Y=\emptyset$, since $R\cup S \subseteq W$.
There is no path of $G[C]$ of length at most $k-2$ between $Y$ and $S$, since $Y\cap W=\emptyset$.
Since $Y\subseteq L_k$, no vertex in $Y$ has $G[C]$-distance less than $k$ from a vertex in $R$; and since
$y\cap W=\emptyset$, it follows 
that every vertex in $Y$ has $G[C]$-distance exactly $k$
from some vertex in $R$. Since $Y\ne \emptyset$, there exists $r\in R$ such that $X\ne \emptyset$, 
where $X$ denotes the set of vertices in $Y$
with $G[C]$-distance $k$ from $r$.
From \ref{bigtwiglemma} applied to $X$ and $G[Y]$,
there exists $x\in X$ that is $d$-equipped in $Y\setminus X$. Let $P$ be a path of $G[C]$ of length $k$ joining $r,x$.
Since no vertex in $Y$ has $G$-distance at most $k$ from a vertex in $\{r_0,b_1\ll b_m\}$, it follows that no vertex of $P$
belongs to $S$. Since $r\in R$, $r$ is adjacent to at least $d$
of $b_1\ll b_m$, say $b_1\ll b_d$. But then adding the parts of the equipment of $x$
to the tree formed by the union of $P$, the paths $Q_1\ll Q_d$ and the edges $rb_1\ll rb_d$, gives a $(k,d)$-binary star and
a $(k,d)$-bristled star.  This proves (1).

\bigskip

Thus we may assume that $\chi(L_k\setminus W)\le c_0$. Let $X=L_0\cup\cdots\cup L_k$. Since $\chi(W)\le (m+1)\tau$, and 
$L_0\cup\cdots\cup L_{k-1}\subseteq W$,
it follows that $\chi(X)\le (m+1)\tau+c_0$.
By \ref{doublecreatures} applied to $G[C]$ and $X$,
there exists $x\in X$ that is $d$-equipped in $C\setminus X$, 
and there exists $x'\in C\setminus X$ adjacent to $x$ such that $x'$ is $d$-equipped in $C\setminus (X\cup N^1(x))$.
Since $x'\notin X$, it follows that $x'\in L_{k+1}$, and $x\in L_k$, and so there is an induced path $T$ of $G[C]$
with ends $x',u$ say, such that $x\in V(T)$, and either
\begin{itemize}
\item $k\ge 2$, and $T$ has length $k-1$ and $u\in S$, or 
\item $T$ has length $k+1$, and $u\in R$.
\end{itemize}
In either case no vertex of $T\setminus u$ is in $R$.

Suppose that the first bullet holds. It follows that no vertex of $T$ different from $u$ is in $S$; and
since $k\ge 2$ and so $S\subseteq X$, it follows that no vertex of the equipment of $x$
belongs to $S$. 
Since $u\in L_2$, it follows that $u$ is nonadjacent to $r_0$, and since $u\in S$ 
we may assume that $u$ is nonadjacent to 
$b_1\ll b_d$ and adjacent to $b_{d+1}$. But then 
the union of the paths $Q_1\ll Q_d$, the edges $r_0b_1\ll r_0b_d$, the edge $r_0b_{d+1}$, 
the edge $b_{d+1}u$, and the path $T$ can be extended
to a $(k,d)$-binary star and to a $(k,d)$-bristled star by adding appropriate parts of the equipment of $x$.

We may therefore assume that the second bullet holds. Let $v$ be the vertex of $T$ adjacent to $u$, and $w$ the other neighbour
of $v$ in $T$ (this exists since $T$ has length $k+1\ge 2$).
Since $x'\in L_{k+1}$, no vertex of $T$ different from $v,w$ belongs to $S$. 
\\
\\
(2) {\em If $S\subseteq X$ 
then the theorem holds.}
\\
\\
Since $S\subseteq X$, it follows that $x'\notin S$. Moreover, if $w\in S$ then 
$T\setminus \{u,v\}$ satisfies the first bullet above and we are done; so we may assume that $w\notin S$. 
We may assume that $u$ is adjacent to $b_1\ll b_{2d}$. 
If $v$ is nonadjacent to at least $d$ of $b_1\ll b_{2d}$, say to $b_1\ll d_d$, then the union of the paths
$Q_1\ll Q_d$, the edges $ub_1\ll ub_d$, the path $T\setminus x'$, and appropriate parts of the equipment of $x$, 
gives a $(k,d)$-binary star and a $(k,d)$-bristled star.
Thus we may assume that $v$ is adjacent to at least $d$ of $b_1\ll b_{2d}$, say to $b_1\ll b_d$. But then 
the union of the paths
$Q_1\ll Q_d$, the edges $vb_1\ll vb_d$, the path $T\setminus u$, and appropriate parts of the equipment of $x'$, 
gives a $(k,d)$-binary star and a $(k,d)$-bristled star. This proves (2).

\bigskip

In view of (2), we may assume that $S\not\subseteq X$, and consequently $k=1$. Choose $s\in S\setminus X$; 
then $s$ is nonadjacent to $r_0$, and we may assume that
$s$ is nonadjacent to $b_1\ll b_d$ and adjacent to $b_{d+1}$, and so the union of the paths $Q_1\ll Q_{d+1}$, the edges
$r_0b_1\ll r_0b_{d+1}$ and the edge $b_{d+1}s$ gives a $(1,d)$-bristled star. It remains to find a $(1,d)$-binary star.
To do so, we need to apply \ref{doublecreatures} more carefully, using the ``forbidden vertices'' feature of \ref{doublecreatures}. 
First we need:
\\
\\
(3) {\em We may assume that $|S\setminus L_1|< mn$.}
\\
\\
If not, then some vertex in $\{b_1\ll b_m\}$ (say $b_{d+1}$) is adjacent to more than $n$ vertices in $S\setminus L_1$, and 
from the definition of $n$, at 
least $d$ of these neighbours are pairwise nonadjacent, say $y_1\ll y_d$. For $1\le i\le d$, $y_i$ is adjacent to fewer than
$2d$ of $\{b_1\ll b_m\}$, and since $m=2d^2$, we may assume that $y_1\ll y_d$ are all nonadjacent to all of $b_1\ll b_d$.
But then the union of the paths $Q_1\ll Q_d$, the edges $r_0b_1\ll r_0b_d, r_0b_{d+1}$, and the edges
$b_{d+1}y_1\ll b_{d+1}y_d$ gives a $(1,d)$-binary star.
This proves (3).

\bigskip

Now let us apply \ref{doublecreatures} again, to $X=L_0\cup L_1$ and $G[C]$, setting $B=S\setminus L_1$ and replacing $d$ by 
$mn+d$. We obtain a pair $x,x'$ as before, where $x$ is $(mn+d)$-equipped (and hence $d$-equipped) in
$C\setminus X$, and $x'\in L_2\setminus S$, adjacent to $x$, 
such that $x'$ is $(mn+d)$-equipped in $S\cup (C\setminus (X\cup N^1(x)))$. In particular, since $x$ has at least
$mn+d$ pairwise nonadjacent neighbours in $L_2$, and at most $mn$ of them belong to $S$, it follows that
$x$ has $d$ pairwise nonadjacent neighbours $y_1\ll y_d\in L_2\setminus S$. Also, $x'\in L_2\setminus S$,
and by the same argument $x'$ has $d$ pairwise nonadjacent neighbours $y_1'\ll y_d'\in C\setminus (X\cup S\cup N^1(x))$.
But now we finish the proof as in (2). More precisely, let $T$ be the path $u\d x\d x'$, where $u\in R$. 
We may assume that $u$ is adjacent to $b_1\ll b_{2d}$. 
If $x$ is nonadjacent to at least $d$ of $b_1\ll b_{2d}$, say to $b_1\ll b_d$, then the union of the paths
$Q_1\ll Q_d$, the edges $ub_i\;(1\le i\le d)$, the edge $ux$, and the edges $xy_1\ll xy_d$
gives a $(1,d)$-binary star.
If $x$ is adjacent to at least $d$ of $b_1\ll b_{2d}$, say to $b_1\ll b_d$, then 
the union of the paths
$Q_1\ll Q_d$, the edges $xb_i\;(1\le i\le d)$, the edge $xx'$, and the edges $x'y_1'\ll x'y_d'$
gives a $(1,d)$-binary star. This proves \ref{freecathedral}.~\bbox

We can extend this result to cathedrals that are not free, as follows.
\begin{thm}\label{cathedral}
For all $\tau\ge 0$ and $k,d\ge 1$ there exist $c,n\ge 0$ with the following property.
Let $G$ be a graph with $\chi^k(G)\le \tau$, and let
$(\mathcal{S}_1\ll\mathcal{S}_n)$ be a cathedral in $G$ of height $d$, dominating a set $C$ with $\chi(C)>c$.
Then $G$ is $(k,d)$-starry.
\end{thm}
\Proof
Choose $c_0,n_0$ such that \ref{freecathedral} holds with $c,n$ replaced by $c_0, n_0$. 
Let $n=dn_0$, and choose $c\ge 2^{n^2}c_0$ such that 
\ref{bigtwiglemma} holds with $a,c$ replaced by $\tau, c$.
We claim that $n,c$ satisfy \ref{cathedral}.  For let 
$G$ be a graph with $\chi^k(G)\le \tau$, and let
$(\mathcal{S}_1\ll\mathcal{S}_n)$ be a cathedral in $G$ of height $d$, dominating a set $C$ with $\chi(C)>c$.
Let $\mathcal{S}_i=(P_i,A_i,B_i)$ for $1\le i\le n$.
We may assume that $C$ is minimal such that $\chi(C)>c$; and for $1\le i\le n$, that every vertex in $B_i$ has a neighbour in $C$
(because any vertex in $B_i$ with no neighbour in $C$ can be removed). 
\\
\\
(1) {\em We may assume that there exists $i$ with $1\le i\le n$, and a vertex $v\in B_i$ that has neighbours in at least $d$
of $A_{i+1}\ll A_n$.}
\\
\\
Suppose not.
For $1\le i\le n$ and each $v\in B_i$, let $J(v)$ be the set of $j$ with $i<j\le n$ such that $v$ has a neighbour in $A_j$.
Thus each $|J(v)|\le d-1$.
For each vertex $u\in C$, choose a neighbour $v_i\in B_i$ of $u$, and let $S(u)$ be the sequence
$(J(v_1)\ll J(v_n))$. There are at most $2^{n^2}$ possibilities for this sequence, and so there exists $C'\subseteq C$
with $\chi(C')\ge \chi(C)2^{-n^2}>c_0$ and a sequence $S$ such that $S(v)=S$ for all $v\in C'$. Let $S=(J_1\ll J_n)$.
Since there exists $v\in B_i$ with $J(v)=J_i$, it follows that $|J_i|<d$ for $1\le i\le n$. Let $H$
be the digraph with vertex set $\{1\ll n\}$ in which $j$ is adjacent from $i$ if $i<j$ and $j\in J_i$. Since $H$ has no
directed cycles and all vertices
have outdegree less than $d$, the graph underlying $H$ has chromatic number at most $d$. Consequently
it has a stable set $I$ of cardinality $n_0$, since $n=dn_0$.
For each $i\in I$, let $B'_i$ be the set of vertices $v\in B_i$ such that $J(v)=J_i$. It follows that $B_i'$ covers $C'$.
But then for each $i\in I$, $(P_i, A_i, B_i')$ is a spire dominating $C'$, and the sequence of these spires is a free cathedral
of height $d$ and length $n_0$ dominating $C'$, and the result follows from \ref{freecathedral}.
This proves (1).

\bigskip

Choose $v,i$ as in (1), where $v \in B_i$ has neighbours in at least $d$
of $A_{i+1}\ll A_n$, say in $A_{j_1}\ll A_{j_d}$. If $v$ has a neighbour in $A_j$ for some $j>i$, then
since $v$ has no neighbour in $V(P_j)\setminus A_j$, there is an induced path $Q_j$ of length $d$
with one end $v$ and with all other vertices in $A_j\cup V(P_j)$. Let $R$ be the union of
the paths $Q_{j_1}\ll Q_{j_d}$. Let $S$ be the set of neighbours of $v$ in $C$; thus $S\ne \emptyset$ since every vertex
of $B_i$ has a neighbour in $C$.
Let $X$ be the set of vertices in $C$
that can be joined to a vertex in $S$ by a path of $G[C]$ of length at most $k-1$. Thus $S\subseteq X$, and $\chi(X)\le \tau$.
By \ref{bigtwiglemma}, there is a vertex $x\in X$ that is $d$-equipped in $C\setminus X$. Choose a path
$T$ of $G[C]$ between $x$ and some $u\in S$, of length $k-1$. Since $x$ has a neighbour in $C$ that does not belong
to $X$, it follows that $T$ has length $k-1$ and no vertex of $T$ except $u$ belongs to $S$. But then the union of the
paths $Q_{j_1}\ll Q_{j_d}$, the edge $vu$, the path $T$, and appropriate parts of the equipment of $x$, gives a
$(k,d)$-binary star and a $(k,d)$-bristled star. This proves \ref{cathedral}.~\bbox

\section{Building a cathedral}

To apply \ref{cathedral} we need to prove that our graph contains an appropriate cathedral.
First we need:

\begin{thm}\label{getspire}
For all $c,d,\tau\ge 0$ there exists $c'\ge 0$ with the following property.
Let $G$ be a graph with $\chi^1(G)\le \tau$ and $\chi(G)>c'$. Then there is a spire
$(P,A,B)$ in $G$ with height $d$, dominating a set $C$ with $\chi(C)>c$.
\end{thm}
\Proof We may assume that $c\ge \tau$, by increasing $c$ if necessary.
Let $c'=2c+d\tau+1$; we claim that $c'$ satisfies the theorem. For let $G$ be as in the theorem, with 
$\chi(G)>c'$. We may assume that $G$ is connected. Choose a vertex $x_0$. Since $\chi(G)>c'$,
there is a component of $G\setminus x_0$ with chromatic number at least $c'$, with vertex set $C_1$ say.
By \ref{getpath} applied to $x_0$ and $C_1$, there is an induced path $x_0 \c x_{d}$ of $G$ 
where $x_1\ll x_{d}\in C_1$, and a subset $C_2$ of $C_1$,
with the following properties:
\begin{itemize}
\item $x_0\ll x_{d}\notin C_2$;
\item $G[C_2]$ is connected;
\item $x_{d}$ has a neighbour in $C_2$, and $x_0\ll x_{d-1}$ have no neighbours in $C_2$; and
\item $\chi(C_2)\ge \chi(C_1)-d\chi^1(G)>2c$.
\end{itemize}
Let $P$ be the path $x_0\c x_d$. For $i\ge 0$, let $L_i$ be the set of vertices in $C_2$
with $G[C_2]$-distance exactly $i$ from $x_d$, and choose $i\ge 0$ such that $\chi(L_i)\ge \chi(C_2)/2> c$.
Since $c\ge \tau$ it follows that $i\ge 2$. Let $A=L_0\cup\cdots\cup L_{i-2}$, and $B=L_{i-1}$; then $(P,A,B)$
is a spire of height $d$ dominating $L_i$. This proves \ref{getspire}.~\bbox

\begin{thm}\label{getcath}
For all integers $d\ge 1$ and $n,c,\tau\ge 0$, there exists $c'\ge 0$ with the following property. 
Let $G$ be a graph such that $\chi^2(G)\le \tau$ and $\chi(G)>c'$.  Then $G$ admits a
cathedral of length $n$ and height $d$ that dominates a set of chromatic number more than $c$.
\end{thm}
\Proof 
For each integer $x\ge 0$, choose $\phi(x)$ such that \ref{getspire} holds with $c,c'$ replaced by $x,\phi(x)$.
Let $c_n=c$, and inductively
for $n-1\ge i\ge 0$ let $c_i=\phi(d\tau+c_{i+1})$. Let $c'=c_0$. 
Now let $G$ be as in the theorem with $\chi(G)>c'$.
For $1\le i\le n$, choose a spire $\mathcal{S}_i=(P_i,A_i,B_i)$ and a set $C_i$ inductively as follows.
Let $C_0=V(G)$. Suppose that for some $i$ with $1\le i\le n$ we have chosen $C_{i-1}$ with $\chi(C_{i-1})>c_{i-1}$.
Since $c_{i-1}=\phi(d\tau+c_{i})$, there is a spire $(P_i,A_i,B_i)$ in $G[C_{i-1}]$, dominating a set $C_i'\subseteq C_{i-1}$
with $\chi(C_i')>d\tau+c_i$. The set of vertices of $C_i'$ with $G$-distance at most two from a vertex of $V(P_i)\setminus A_i$
has chromatic number at most $d\tau$, so there exists $C_i\subseteq C_i'$ with $\chi(C_i)>c_i$, such that 
every path in $G$ between $C_i$ and a vertex of $V(P_i)\setminus A_i$ has length at least three.
This completes the inductive definition of $\mathcal{S}_i=(P_i,A_i,B_i)$ for $1\le i\le n$.
We see that:
\begin{itemize}
\item the spires $\mathcal{S}_i\;(1\le i\le n)$ are pairwise vertex-disjoint, and each of them dominates $C_n$; and
\item for $1\le i<j\le n$, if $u\in V(\mathcal{S}_i)$ is adjacent to $v\in V(\mathcal{S}_j)$ then $u\in B_i$.
\end{itemize}
For $1\le i\le n$, let $B_i'$ be the set of vertices in $B_i$ that have neighbours in $C_n$. Since
for $1\le j\le n$ every path in $G$ between $C_n$ and a vertex of $V(P_j)\setminus A_j$ has length at least three,
it follows that for $1\le i<j\le n$, no vertex in $B_i'$ has a neighbour in $V(P_j)\setminus A_j$.
Let $\mathcal{S}_i'=(P_i, A_i, B_i')$ for $1\le i\le n$; then 
$(\mathcal{S}_1'\ll\mathcal{S}_n')$ is a cathedral in $G$ of height $d$, dominating $C_n$, and $\chi(C_n)>c$.
This proves \ref{getcath}.~\bbox

\bigskip

Combining \ref{cathedral} and \ref{getcath} we obtain:
\begin{thm}\label{cath2}
For all $k,d\ge 1$ and $\tau\ge 0$ there exists $c$ such that, if $G$ is a graph with $\chi^2(G)\le \tau$ and $\chi^k(G)\le \tau$,
and $\chi(G)>c$, then
$G$ is $(k,d)$-starry.
\end{thm}

\section{The case $k=1$}

It remains to prove \ref{cath2} with the hypothesis that $\chi^2(G)\le \tau$ omitted.
(We may assume that it is not implied by the hypothesis $\chi^k(G)\le \tau$, and so we only need to handle the case $k=1$.)
The proof of \ref{getcath} no longer works, since in the notation of \ref{getcath}
we have no way to stop
vertices in $B_i'$ having neighbours in later paths $P_j$.
The content of this section is our workaround.

If $H$ is a $d$-superstar, we call its 
vertex of degree $d$ its {\em root}.
From repeated application of \ref{treesplit} with $r=2$ and $k=d$ we deduce:

\begin{thm}\label{starsplit}
For all $c,d,\tau\ge 0$,  there exists $c'\ge 0$ with the following property. Let $G$ be a graph with $\chi^1(G)\le \tau$, and let
$z\in V(G)$ such that
$\chi(N^2_G(z))>c'$. Then there is a $d$-superstar
$H$
with root $z$, and a subset $W$ of $V(G)$, such that
\begin{itemize}
\item $z\in W$, and no other vertex of $H$ belongs to $W$ or has a neighbour in $W\setminus \{z\}$; and
\item $\chi(N^2_{G'}(z))>c$, where $G'$ is the graph $G[W]$.
\end{itemize}
\end{thm}
Let us say a {\em $d$-band} in $G$ is a triple $(H,z,B)$, where 
\begin{itemize}
\item $H$ is a $d$-superstar in $G$ with root $z$;
\item $B\cap V(H)=\emptyset$, and $z$ is adjacent to every vertex in $B$; and
\item no vertex in $H\setminus \{z\}$ has a neighbour in $B$.
\end{itemize}
(Thus, a band is like a spire, but with the path replaced by a subdivided star, and the set $A$ is just a single vertex.)
If $\mathcal{S} = (H,z,B)$ is a band, we write $V(\mathcal{S})=V(H)\cup B$. A $d$-band $(H,z,B)$ {\em dominates} $C$ if 
$V(H)\cap C=\emptyset$, and $B$ covers $C$, and there is no edge between $V(H)$ and $C$.
We deduce:

\begin{thm}\label{k=1proof}
For all $d\ge 1$ and $\tau\ge 0$ there exists $c\ge 0$ with the following property. Let $G$ be a graph with
$\chi^1(G)\le \tau$  and $\chi(G)>c$.
Then $G$ is $(1,d)$-starry.
\end{thm}
\Proof
Choose $n_0,c'$ such that \ref{freecathedral} is satisfied setting $k=1$, $n=n_0$ and $c=c'$.
Let $n=(2d+1)n_0$.
For each $x\ge 0$, choose $\phi(x)$ such that,
if $G$ is a graph with $\chi^2(G)\le x$
and $\chi(G)>\phi(x)$, then
$G$ is $(1,d)$-starry. (This is possible by \ref{cath2}.)
For each $x\ge 0$, choose $\psi(x)$ such that \ref{starsplit} holds with $c,c'$ replaced by $x,\psi(x)$. 

Let $c_n=\max(c'2^{n^2},d\tau)$, and
for $n-1\ge i\ge 0$ let $c_i=\phi(\psi(c_{i+1}))$. Let $c=c_0$.
Now let $G$ be as in \ref{k=1proof} with $\chi(G)>c$.
For $1\le i\le n$, choose a $d$-band $\mathcal{S}_i=(H_i,z_i,B_i)$ and a set $C_i$ inductively as follows.
Let $C_0=V(G)$. Suppose that for some $i$ with $1\le i\le n$ we have chosen $C_{i-1}$ with $\chi(C_{i-1})>c_{i-1}$.
Since $c_{i-1}\ge \phi(\psi(c_i))$, we may assume that $\chi^2(G[C_{i-1}])>\psi(c_{i})$, for otherwise the result follows from
\ref{cath2}. By \ref{starsplit}, 
there is a $d$-band $(H_i,z_i,B_i)$ in $G[C_{i-1}]$, dominating a set $C_i'\subseteq C_{i-1}$
with $\chi(C_i')>c_i$. 
This completes the inductive definition of $\mathcal{S}_i=(H_i,z_i,B_i)$ for $1\le i\le n$.

We see that:
\begin{itemize}
\item the $d$-bands $\mathcal{S}_i\;(1\le i\le n)$ are pairwise vertex-disjoint, and each of them dominates $C_n$; and
\item for $1\le i<j\le n$, if $u\in V(\mathcal{S}_i)$ is adjacent to $v\in V(\mathcal{S}_j)$ then $u\in B_i$.
\end{itemize}
We may assume that $G[C_n]$ is connected, and every vertex in $B_i$ has a neighbour in $C_n$, for $1\le i\le n$.
\\
\\
(1) {\em We may assume that there exists $i$ with $1\le i\le n$, and a vertex $v\in B_i$ that has neighbours in at least $d$
of $V(H_{i+1})\ll V(H_n)$.}
\\
\\
Suppose not; then as in the proof of step (1) of \ref{cathedral}, there exist $I\subseteq \{1\ll n\}$ with 
$|I|=n_0$, and $B_i'\subseteq B_i$ for each $i\in I$, and a subset $C'\subseteq C_n$ with $\chi(C')\ge \chi(C_n)2^{-n^2}>c'$,
such that $B_i'$ covers $C'$ for each $i\in I$, and for $i,j\in I$ with $i<j$, there are no edges between $B_i'$
and $V(H_j)$. For each $i\in I$, let $P_i$ be a path of $H_i$ of length $d$ with one end $z_i$; then
$(P_i,\{z_i\},B_i')$ is a spire dominating $C'$, and the result follows from \ref{freecathedral}.
This proves (1).

\bigskip

By (1), and by renumbering, we may assume that there exists $v \in B_1$ with neighbours in 
$A_{2}\ll A_{d+1}$. 
(We no longer need the other $d$-bands, so this renumbering is legitimate.)
For $2\le i\le d+1$ let $y_i\in V(H_i)$ be adjacent to $v$. Then the union of $H_1$, the edge
$z_1v$, and the edges $vy_2\ll vy_{d+1}$ forms a $(1,d)$-binary star. Since $v$ has a neighbour in $C_n$,
and $\chi(C_n)>d\tau$, \ref{getpath} implies that there is an induced path $P$ of length $d$ with one end $v$
and all other vertices in $C_n$. But then the union of $H_1$, the edges $z_1v$ and $vy_2$, and the path $P$,
forms a $(1,d)$-bristled star. This proves \ref{k=1proof}.~\bbox

Let us put these pieces together to deduce \ref{mainthm}.

\bigskip

\noindent{\bf Proof of \ref{mainthm}.\ \ }
Let $k,d\ge 1$, and $\kappa\ge 0$; we need to prove that for all $\kappa\ge 0$, there exists
$c$ such that 
every graph $G$ with $\omega(G)\le \kappa$ and $\chi(c)>c$ is $(k,d)$-starry.
We prove this by induction on $\kappa$.
Let $\tau_1$ be such that
every graph $G$ with $\omega(G)\le \kappa-1$ and $\chi(G)>\tau_1$ is $(k,d)$-starry.
By \ref{bigball} with $r=k$, there exists $\tau\ge \tau_1$ such that 
every graph $G$ with $\chi^1(G)\le \tau_1$ and $\chi^k(G)>\tau$ is $(k,d)$-starry.
Choose $c_1$ such that \ref{cath2} holds with $c$ replaced by $c_1$.
Choose $c_2$ such that \ref{k=1proof} holds with $\tau,c$ replaced by $\tau_1,c_2$.
Let $c=\max(c_1,c_2)$, and let $G$ be a graph with $\chi(G)>c$ and $\omega(G)\le \kappa$. We may assume that
$\chi^1(G)\le \tau_1$, for otherwise the result follows from the induction on $\omega(G)$.
We may assume that $\chi(G^k)\le \tau$ for otherwise the result follows from \ref{bigball}.
By \ref{cath2} we may assume that $\chi^2(G)> \tau$ and so $k=1$ because $\chi(G^k)\le \tau$. But then the result
follows from \ref{k=1proof}. This proves \ref{mainthm}.~\bbox

\section{Two counterexamples}

As we said in the beginning, our main tools are the lemmas of section 2. We proved there in particular that if $G$ is a graph 
of very large chromatic number and with $\chi^1(G)$ bounded, and $v$ is a vertex such that $\chi(G\setminus v)<\chi(G)$, 
then there is an induced $d$-star in $G$ with centre $v$, and there is an induced $(d+1)$-edge path in $G$ with second vertex $v$.
What other trees containing $v$ must be present? 
If we could find more, then 
the methods of this paper might allow us to prove the Gy\'arf\'as-Sumner conjecture for more types of trees.
But there are not many more. For instance, Sophie Spirkl and the third author showed that there need not be an induced 
five-vertex path in $G$ with middle vertex $v$. Here is the example.

Choose a large integer $k$, and take a minimal triangle-free graph with chromatic number more than $k$. Let $I$
be the set of neighbours of some vertex $u$, and delete $u$. This produces a graph $H$ say, and a stable subset $I$ of $V(H)$,
such that $H$ is triangle-free, and $k$-colourable, but in every $k$-colouring all $k$ colours occur in $I$. 

For each subset $S\subseteq I$ with $|S|=k-1$, make a gadget $B^S$ as follows. Take some enumeration $\{s_1\ll s_{k-1}\}$
of $S$, take $2k$ new vertices
$a_1^S,b_1^S,a_2^S,b_2^S\ll a_{k}^S, b_{k}^S$, and for $1\le i\le k$ and $i\le j\le k-1$ 
make $a_i^S$ and $b_i^S$ both adjacent to $s_j$. 
Also for $1\le i<j\le k-1$, make $a_i^S$ adjacent to $b_j^S$ and $b_i^S$ adjacent to $a_j^S$. It is easy to see
that every $k$-colouring of $S$ can be extended to a $k$-colouring of the new vertices. Now add one more new vertex
$v^{S}$ adjacent to all of $a_1^S,b_1^S\ll a_{k}^S, b_{k}^S$. This defines $B^S$ (that is, the set
$\{a_1^S,b_1^S\ll a_{k}^S, b_{k}^S, v^S\}$, and the new edges incident with these new vertices). Let $G_1$ be the graph obtained from
$H$ by adding $B^S$ for every choice of $S$. Now $G_1$ is not $k$-colourable; because in a $k$-colouring of $H$, some choice of $S$
is coloured with all different colours, and this cannot be extended to a $k$-colouring of $G_1$. 
Starting with $H$, let us add the gadgets $B^S$ one by one until the chromatic number increases to $k+1$, and then stop; let
$G$ be the graph just constructed and $B^S$ the final gadget added. Let $v=v^S$; then $\chi(G\setminus v)<\chi(G)$, and
there is no five-vertex induced path of $G$ with middle vertex $v$.

One might also hope that the parts of the equipment of a vertex can be unified; 
say a vertex $v$ is ``properly $d$-equipped'' in $Y$ if there is an induced path of length $d$ with first vertex $v$
and all other vertices in $Y$, and $d$ pairwise nonadjacent neighbours of $v$, all in $Y$, and such that none of them has
any neighbours in $P\setminus v$. One might hope that \ref{bigstabletwiglemma} could be strengthened correspondingly.
But this is false, even for $d=2$ and for triangle-free graphs, 
as another counterexample (also due to Sophie Spirkl and the third author) shows. Take $H$ and $I$ as before;
and for each $S\subseteq I$ with $|S|=k$, let $B^S$ be a gadget defined as follows. Let $S=\{s_1\ll s_k\}$, and take
$k+1$ new vertices $a_1^S\ll a_{k}^S$ and $v^S$, and for all distinct $i,j$ with $1\le i,j\le k$, make $a_i^S$ adjacent to 
$s_j$. Also make $v^S$ adjacent to $a_1^S\ll a_{k}^S$. Adding all these gadgets increases the chromatic number, so, as before,
add them one at a time until the chromatic number increases, let $B^S$ be the last one, and let $v=v^S$. Then $v$ is
not properly $2$-equipped in $V(G)\setminus \{v\}$.

\end{document}